\documentclass[11pt,twoside]{article}
\usepackage[cp1251]{inputenc}
\usepackage[russian]{babel}
\usepackage{amssymb,amsmath}
\usepackage[dvips]{graphicx}

\newcommand{\Title}[1]{\begin{center}\baselineskip=6.0mm{\Large\textbf{#1}}\end{center}}
\newcommand{\Author}[1]{\centerline{\large\textbf{#1}}}

\newcommand{\Email}[1]{\small\texttt{#1}}
\newcommand{\lit}[3]{\vspace*{0.7mm}\par\noindent\makebox[5.2mm][r]{#1.~}\parbox[t]{162.8mm}{{\textit{#2}}~{#3}}\hspace*{-1.6mm}}

\begin{document}

\pagestyle{myheadings}

\sloppy
\Title{Values of the $t$-invariant for small Seifert manifolds\footnote{The work is supported by
RFBR (grant ~06-01-72014).}}

\Author{Mikhail Ovchinnikov}

Chelyabinsk State University, Chelyabinsk;

\Email{ovch@csu.ru}

\bigskip
\bigskip

{\footnotesize

The $t$-invariant can be considered as the Turaev-Viro invariant of order 5 computed for integer colors only.

We compute all values of the $t$-invariant for Seifert manifolds with base sphere and three singular fibers. As a result we show that
the manifolds parameters modulo five define the value of the $t$-invariant. Partially we show that there are 12  distinct values of the $t$-invariant for these manifolds. Some examples show that the $t$-invariant for these manifolds is not defined by the first homology group.
}

\vspace{3mm}

{\bf Introduction}

\vspace{1mm}

A polyhedron $P$ in a compact 3-manifold $M$ is called {\em a spine} of the manifold if its complement in $M$ is homeomorphic to the product of the manifold boundary and halfinterval or  to open 3-ball. A two-dimensional polyhedron is called {\em simple polyhedron} if each its point possesses a neighbourhood embeddable in the space
$$ (R\times \{0\}\times R_+ )\cup (R^2\times\{0\} \cup(\{0\}\times R \times R_-)$$

If the point is mapped to the point $(0,0,0)$ by such embedding then it is called {\em a vertice of the polyhedron}.
A spine is called  {\em a simple spine} if it is a simple polyhedron  ~([1]).

{\em $t$-invariant} of a compact 3-manifold $M$ is defined by the formula
$$ t(M)=\sum_{Q\subset P}{(-\varepsilon)^{-v(Q)}\varepsilon^{\chi(Q)}},$$
where $Q$ denotes a simple subpolyhedron, possibly empty, in a simple spine  $P$, $v(Q)$ is the number of vertices of the subpolyhedron , $\chi(Q)$~ is the Euler characteristic of the subpolyhedron, $\varepsilon$ is a root of the equation $\varepsilon^2=\varepsilon+1$.

The  $t$-invariant is a topological invariant of a manifold because its value does not depend on the choice of a simple spine  $P$ of the manifold  ~([2]).

It is known that  the $t$-invariant can be considered as a slight modification of the Turaev-Viro invariant for a 5-th degree root of unity. The modification is a restriction of the construction of the Turaev-Viro invariant for integer colors~([2]).

\vspace{3mm}

\textbf{Theorem}
{\em
Let $M = (b, (O_1,0), (\alpha_1,\beta_1), (\alpha_2,\beta_2), (\alpha_3,\beta_3))$ be a small Seifert manifold.

If $\alpha_3 = 1  \,\rm{mod} 5$, then $t(M)=t(L_{p,q})= \left\{ \begin
{array}{rl}
0,            & \mbox{if } p=0  \,\rm{mod} 5, q=\pm 1  \,\rm{mod} 5 \\
1,            & \mbox{if } p=\pm 1  \,\rm{mod} 5 \\
\varepsilon,  & \mbox{if } p=\pm 2  \,\rm{mod} 5 \\
\varepsilon+2,& \mbox{if } p=0  \,\rm{mod} 5, q=\pm 2  \,\rm{mod} 5
\end{array} \right.,$

where $p=\alpha_1\beta_2 - \alpha_2\beta_1 + \alpha_1\alpha_2(\beta_3-b), q=\alpha_1\lambda+\beta_1\mu, \lambda$ and $\mu$ are any numbers satisfying $\alpha_2\lambda + \beta_2\mu = 1$.

If $\alpha_3 = 0  \,\rm{mod} 5, \beta_3 = \pm 1   \,\rm{mod} 5$, then $t(M)=t(L_{\alpha_1,\beta_1})t(L_{\alpha_2,\beta_2})$.

If $\alpha_i = \pm 2   \,\rm{mod} 5, i=1,2,3$, then
$t(M)= \left\{ \begin
{array}{rl}
 \varepsilon+2,& \mbox{if } k=0  \,\rm{mod} 5 \\
2\varepsilon+3,& \mbox{if } k=\pm 1  \,\rm{mod} 5 \\
 \varepsilon+3,& \mbox{if } k=\pm 2  \,\rm{mod} 5
\end{array} \right.,$
where $k = b - 2(\beta_1+\beta_2+\beta_3)  \,\rm{mod} 5 $.

If $n \ge 1$ singular fibers have parameters $0  \,\rm{mod} 5, \pm 2  \,\rm{mod} 5 $, and the first parameter of any other singular fiber equals $\pm 2  \,\rm{mod} 5$, then

\centerline{
$t(M)= \left\{ \begin
{array}{rl}
2 - \varepsilon,& \mbox{if } n=1 \\
3 - \varepsilon,& \mbox{if } n=2 \\
5              ,& \mbox{if } n=3
\end{array} \right.$
}

}

\vspace{2mm}

A scheme of the proof. Firstly we describe two constructions which we use essentially in our proof. The first one is a generalized  $t$-invariant for a simple polyhedron with boundary with simple subgraph in the boundary chosen. The second one is a method to build spines out of copies of some elementary polyhedra. The constructions imply  a presentation of values of the  $t$-invariant as  products of some matrices. Studying these products we obtain that repasting of a solid torus in a 3-manifold gives at most 12 distinct values of the $t$-invariant and there are some equations between these 12 values.  Also we use equivalences for different presentations of the same Seifert manifolds. That occurs to be enough to compute the necessary values of the $t$-invariant by short calculations.

\vspace{4mm}

{\bf A generalization of the $t$-invariant for polyhedra with boundary with a simple subgraph chosen}

\vspace{1mm}

Definition. A polyhedron is called  {\em a simple polyhedron with boundary} if each its point has a neighbourhood embeddable in the space
$$ (R_+\times \{0\}\times R_- )\cup (R_+\times R\times\{0\} \cup (\{1\}\times R \times R_+)$$

Points which are mapped by such embeddings to  $ (\{0\}\times \{0\}\times R_- )\cup (\{0\}\times R\times\{0\})$  form  {\em the boundary} of the polyhedron.

A graph is called {\em a simple graph} if all its vertices are of valence three.

{\em A generalized $t$-invariant} of simple polyhedron $P$ with boundary with a simple subgraph  $G$ chosen in the boundary is defined by the formula
$$ t(P,G)=\sum_{Q\subset P}{(-\varepsilon)^{-v(Q)}\varepsilon^{\chi(\rm{Int}Q)+\frac{1}{2}\chi(\partial{Q})}},$$
\noindent
where $Q$ denotes a simple polyhedron (probably empty) with boundary  $\partial{Q} = G$,  $v(Q)$ is the number of vertices of the subpolyhedron, $\chi(Q)$~ is the Euler characteristic of the subpolyhedron, $\varepsilon$ is a root of the equation $\varepsilon^2=\varepsilon+1$.

{\bf Remark.} The generalized $t$-invariant is defined here for a polyhedron but not for a manifold. It is sufficient here that for the case of empty boundary the generalized $t$-invariant coincides with the original $t$-invariant.

\vspace{1mm}

The main property of the generalized $t$-invariant used in this work is the equation
$$ t(P \cup Q, G \sqcup H) = \sum_{K \subset P \cap Q}{t(P,G \sqcup  K)t(Q,K \sqcup H)}$$
\noindent
where the intersection  $P \cap Q = \partial{P} \cap \partial{Q}$
of simple polyhedra with boundary  $P$ and $Q$ consists of several connected components of the boundaries of these polyhedra,
$G$, $H$, $K$ are simple subgraphs (possibly empty) in the  graphs $\partial (P \slash Q)$, $\partial (Q \slash  P)$, $P \cap Q$  respectively.

This property follows easily from the definition of the generalized $t$-invariant, and it is a partial case of the important property of all Turaev-Viro invariants~([1]).

\vspace{3mm}

{\bf Construction of small Seifert manifolds simple spines out of polyhedra of three types}

\vspace{1mm}

We use the method of construction of simple spines for Waldhausen manifolds described in ~[3].
It occurs that simple spines of small Seifert manifolds can be constructed out of {\em elementary polyhedra } of three types only. Describe these polyhedra.

The description is constructive. Take three copies of the projective plane with one, two and three holes respectively. Choose for each hole a proper nonsplitting arc in the surfaces such that in each surface these arcs intersect pairwise in a single point. Also we need that in the thriceholed surface the arcs do not bound a triangle. These conditions define arcs unambigously.

Attach six halfdisks by their diameters to the surfaces along the arcs. Denote by $E$, $J$, $T$ the polyhedra obtained. Each hole and the corresponding halfdisk give a theta-curve in the boundary of the polyhedron. So these polyhedra have boundary consisting of one, two and three theta-curves respectively.
Now we describe how to built a simple spine of small Seifert manifold with parameters $b, \alpha_1,\beta_1, \alpha_2,\beta_2, \alpha_3,\beta_3$ out of copies of polyhedra $E$, $J$, $T$ by the method described in  ~[3].

Enumerate the edges of each theta-curve by numbers 1, 2, 3 by the following. Number 3 we assign to the boundary edges belonging to the halfdisks. Other two edges of each theta-curve get numbers 1 and 2. Two possible choices of the enumeration in the polyhedron $E$ are equivalent due to the polyhedron symmetry. The twiceholed projective surface is splitted by the arcs into two parts each incident to both boundaries. Assign number 1  to the edges of one part and number 2 to edges of another part.  Two possible choices are equivalent due to symmetry of $J$. The thriceholed projective plane in the polyhedron $T$ is splitted by arcs into four parts. One part is incident to all three holes and each other part is incident to one hole. Assign number 1 to the edges belonging to the first part and assign number 2 to the edges belonging to other parts.

An identifying of two theta-curves with enumerated edges is defined by a permutation of order three and by a correspondence of the vertices. However the symmetries of the polyhedra $E$ and $J$ and the scheme of the pastings imply that the result does not depend on the correspondence of vertices. Also if the permutation is transposition then it is not necessary to specify a direction of the identification.

The spine is obtained from one copy of the polyhedron $T$, three copies of the polyhedron $E$ and few copies of the polyhedron $J$. Each partial pasting is defined by permutation (23) either (13). Now we describe how the pasting is defined by parameters of the manifold. Among equivalent parameters  take parameters satisfying to the  conditions $b = -1, 0<\beta_1 < \alpha_1, 0 < \beta_2 < \alpha_2, 0 < \alpha_3, 0< \beta_3 $. A possibility to provide that follows easily from the theory of Seifert manifolds.  A pair of a singular fiber parameters determines a sequence of pastings for one polyhedron $E$, a few polyhedra $J$ and one theta-curve of the polyhedron $T$. Represent this pair of parameters as the product of few matrices $A$ and $B$ and vector (1,1). The matrix $A$ saves the first coordinate and replaces the second coordinate with the sum of coordinates. The matrix $B$ replaces the first coordinate with the sum of coordinates and saves the second coordinate. Replacing letter $A$ with permutation (23) and letter $B$ with permutation (13) gives the needed pastings. Here the matrix which is multiplyed on vector (1,1) gives the pasting for the polyhedron $E$ and the matrix on other end of the product gives the pasting for theta-curve of the polyhedron $T$.

Associate also two polyhedra to the permutations (23) and (13). These polyhedra are direct product of theta-curve and segment. Each 2-component of such polyhedron is incident to two edges of different theta-curves. In the first polyhedron assign number 1 to edges of one 2-component and numbers 2 and 3 to edges of each other 2-component. In the second polyhedron assign number 2 to edges of one 2-component and numbers 1 and 3 to edges of each other 2-component. These two polyhedra are also called {\em elementary polyhedra}. Realize each gluing by the corresponding polyhedron. Now the construction scheme consists of elementary polyhedra only because now all gluings are the order three identity permutation.

The correctness of the construction is shown in [3].

\vspace{3mm}

{\bf Matrices corresponding to the elementary polyhedra}

\vspace{1mm}

Fix the following order in the set of all five simple subgraphs of theta-curve with numbered edges: empty subgraph, theta-curve without the 1-st edge, theta-curve without the 2-nd edge, theta-curve without the 3-rd edge, theta-curve itself.

Then the values of generalized $t$-invariant for polyhedra $E$ and $J$ form a vector and a matrix unambigously in a natural way.

$ \Phi_E = \left(\begin{matrix} 1 \cr 0 \cr 0 \cr 1 \cr  \varepsilon^{{1}\over{2}} \end{matrix}\right),$
$ \Phi_J = \left(\begin{matrix} 1&0&0&0&0\cr 0&1&0&0&0\cr 0&0&1&0&0\cr
0&0&0&   \varepsilon^{-1}           & \varepsilon^{-{{1}\over{2}}} \cr
0&0&0& \varepsilon^{-{{1}\over{2}}} &-\varepsilon^{-1}  \end{matrix}\right),$

The values of generalized $t$-invariant for polyhedron $T$ form the following cubic matrix $\Phi_T$ presented here as five square matrices.

$ \left(\begin{matrix}
 1&0&0&0&0 \cr
 0&1&0&0&0\cr
 0&0&0&0&0\cr
 0&0&0&0&0\cr
 0&0&0&0&0                         \end{matrix}\right),$
$ \left(\begin{matrix}
 0&1&0&0&0 \cr
 1&1&0&0&0\cr
 0&0&0&0&0\cr
 0&0&0&0&0\cr
 0&0&0&0&0                         \end{matrix}\right),$
$ \left(\begin{matrix} 0&0&0&0&0\cr  0&0&0&0&0\cr
 0&0&0&\!\!\!\! \varepsilon^{-1}\!\!\! &0\cr
 0&0& \!\!\!\! \varepsilon^{-1} \!\!\!\! &0&0 \cr
 0&0&0&0&\!\!\!\! \varepsilon^{-1}\!\!\!\!   \end{matrix}\right),$
$ \left(\begin{matrix} 0&0&0&0&0\cr  0&0&0&0&0\cr
 0&0& \!\!\!\!\! \varepsilon^{-1} \!\!\!\!\! &0&0\cr
 0&0&0& \!\!\!\!\!  \varepsilon^{-2} \!\!\!\!\!  & \!\!\!\!\!  \varepsilon^{-{{3}\over{2}}} \!\! \cr
 0&0&0& \!\!\!\!\!  \varepsilon^{-{{3}\over{2}}} \!\!\!\!\!  & \!\!\!\!\!  -\varepsilon^{-2} \!\! \end{matrix}\right)$
$ \left(\begin{matrix}   0&0&0&0&0\cr 0&0&0&0&0\cr
 0&0&0&0& \!\!\!\! \varepsilon^{-1} \!\!\!\! \cr
 0&0&0& \!\!\! \varepsilon^{-{{3}\over{2}}} \!\!& \!\!\!\!\! -\varepsilon^{-2} \!\!\!\!\! \cr
 0&0& \!\!\! \varepsilon^{-1} \!\!\! & \!\!\! -\varepsilon^{-2} \!\! &  \!\!\!\!\! -\varepsilon^{-{{7}\over{2}}} \!\!\!\!\!
\end{matrix}\right)$

The matrices realizing gluings (23) and (13) occur to be matrices of permutations.

$  \Phi_{(23)} = \left(\begin{matrix}
 1&0&0&0&0\cr
 0&1&0&0&0 \cr
 0&0&0&1&0\cr
 0&0&1&0&0\cr
 0&0&0&0&1 \end{matrix}\right),$
$  \Phi_{(13)} = \left(\begin{matrix}
 1&0&0&0&0\cr
 0&0&0&1&0\cr
 0&0&1&0&0\cr
 0&1&0&0&0\cr
 0&0&0&0&1 \end{matrix}\right)$

{\bf Remark.}
Note that in general case multidimensional matrices appear in such construction  somewhat artificially. At the first step we equip values of the generalized $t$-invariant with corresponding subgraphs of the polyhedron boundary as graphical indexes.  At the second step we group  boundary graphs into disjoint subsets and therefore split the subgraphs into several corresponding subgraphs. At the third step we order subgraphs for each set of boundary components and therefore turn values of indexes from graphs into positive integers. At the fourth step we order  the set of the sets of boundary components and therefore order indexes. Only after that we get a well-defined multidimensional matrix.

A boundary pasting of polyhedra induces the splitting of the set of all indexes of all corresponding matrices into pairs of corresponding indexes. A polyhedron obtained by pasting two polyhedra has corresponding matrix which equals to the product of matrices corresponding to the polyhedra. The product is taken over indexes corresponding to the identified boundary components.

To get the resulting multidimensional matrix we again must  order its indexes.

The symmetry of all square matrices corresponding to the elementary polyhedra and the same symmetry of the cubic matrix $T$ allow not to specify a correspondence between the indexes and the boundaries of the elementary polyhedra.

Of course we can use tensors instead matrices. For that we must choose which boundary of two ones identified must correspond to upper index and then the other one must correspond to low index. But this additional structure seems not to be useful in our construction. The languages of matrices and of tensors are not equivalent. We can find examples of various correspondings of matrix rows and columns to up and low indexes.

The matrix language occured here to be more convenient and clear.

\vspace{3mm}

{\bf Results of matrices products}

Multiply each square matrix and the vector corresponding to the polyhedron $E$ and repeat the multiplication for all vectors which we obtain up to moment when we can't get new vector. We  get 12 vectors.

$\left(\begin{matrix} 1\cr \ 0\ \cr 0\cr \varepsilon\cr 0\end{matrix}\right) \! , \! $
$\left(\begin{matrix} 1\cr \ 0\ \cr \varepsilon\cr 0\cr 0\end{matrix}\right) \! , \! $
$\left(\begin{matrix} 1\cr \varepsilon\cr \ 0\ \cr 0\cr 0\end{matrix}\right) \! , \! $
$\left(\begin{matrix} 1\cr 0\cr 1\cr 0\cr \ \varepsilon^{{{1}\over{2}}} \ \end{matrix}\right) \! , \! $
$\left(\begin{matrix} 1\cr 1\cr 0\cr 0\cr \ \varepsilon^{{{1}\over{2}}} \ \end{matrix}\right) \! , \! $
$\left(\begin{matrix} 1\cr 0\cr 0\cr 1\cr \ \varepsilon^{{{1}\over{2}}} \ \end{matrix}\right) \! , \! $
$\left(\begin{matrix} 1\cr 0\cr 1\cr 1\cr \!\! -\varepsilon^{-{{1}\over{2}}} \!\! \end{matrix}\right) \! , \! $
$\left(\begin{matrix} 1\cr 1\cr 0\cr 1\cr \!\! -\varepsilon^{-{{1}\over{2}}} \!\! \end{matrix}\right) \! , \! $
$\left(\begin{matrix} 1\cr 1\cr 1\cr 0\cr \!\! -\varepsilon^{-{{1}\over{2}}} \!\! \end{matrix}\right) \! , \! $
$\left(\begin{matrix} 1\cr 1\cr \!\! -\varepsilon^{-1} \!\! \cr 1\cr \varepsilon^{-{{3}\over{2}}}\end{matrix}\right) \! , \! $
$\left(\begin{matrix} 1\cr \!\! -\varepsilon^{-1} \!\! \cr 1\cr 1\cr \varepsilon^{-{{3}\over{2}}}\end{matrix}\right) \! , \! $
$\left(\begin{matrix} 1\cr 1\cr 1\cr  \!\!-\varepsilon^{-1} \!\! \cr \varepsilon^{-{{3}\over{2}}}\end{matrix}\right)$

The same result we get if we find all vectors corresponding to singular fibers by the algorithm of constructing spines for small Seifert manifolds.
For this we must start with vectors $\Phi_{(23)}\Phi_E$, $\Phi_{(13)}\Phi_E$ and apply only matrices $\Phi_{(23)}\Phi_J$ and $\Phi_{(23)}\Phi_J$. These 12 points are vertices of an icosahedron.  It is easy to verify that a sort of the distance between distinct vertices is defined by the determinant of the matrix consisting of the vertices parameters $\alpha_1, \beta_1, \alpha_2,\beta_2$. Namely the distance is maximal if the determinant is 0, minimal if the determinant is $\pm 2$, and medial if the determinant is $\pm 1$.

As a result we get that vector corresponding to a singular fiber is defined by the parameters taken modulo five. These 12 vectors correspond to the parameters modulo five $\pm(1,1)$, $\pm(1,0)$, $\pm(0,1)$, $\pm(1,2)$, $\pm(2,1)$, $\pm(1,-1)$, $\pm(1,-2)$, $\pm(2,-1)$, $\pm(2,-2)$, $\pm(2,0)$, $\pm(0,2)$, $\pm(2,2)$ respectively. Already now the finiteness of number of all values of $t$-invariant for small Seifert manifolds is clear. An upper bound for number of the products is $12^3$. The instant number is smaller already due to the symmetries of the presentation. In fact the number of distinct values equals 12.

So the product of the cubic matrix and of three vectors corresponding to parameters $(\lambda_i,\mu_i) \in Z_5+Z_5\backslash \{(0,0)\},i=1,2,3$ is the value of the $t$-invariant for each small Seifert manifold admitting presentation with parameters $b=-1 \bmod{5}, \alpha_i=\lambda_i \bmod{5}, \beta_i = \mu_i \bmod{5}, i=1,2,3$.
To find values of $t$-invariant it is sufficient to find product for the cubic matrix and each triple of vectors taken from the 12 vectors.  It is not hard to calculate all this products and verify the statement of the theorem. Indeed we will decrease radically necessary calculations.

\vspace{3mm}

{\bf A proof of the theorem}

Recall some equivalences of presentations of Seifert manifolds restricted for small Seifert manifolds.
Two small Seifert manifolds $(b, (O_1,0), (\alpha_1,\beta_1), (\alpha_2,\beta_2), (\alpha_3,\beta_3))$ and
$(b', (O_1,0), (\alpha_1',\beta_1'), (\alpha_2',\beta_2'), (\alpha_3',\beta_3'))$ are homeomorphic iff $\alpha_i = \pm \alpha_i', i=1,2,3$ and
$ b+\sum_{i=1}^3{\frac{\beta_i}{\alpha_i}} = b' + \sum_{i=1}^3{\frac{\beta_i'}{\alpha_i'} }$.

That immediately implies that if some $\alpha_i = 0 \bmod{5}$ then the same value of the $t$-invariant we get for a manifold with
$(\alpha_j',\beta_j') = (0, \beta_j \bmod{5})$, if $\alpha_j = 0 \bmod{5}$, and  with $(\alpha_j',\beta_j') = (\alpha_j \bmod{5}, 0) $, if $\alpha_j \ne 0 \bmod{5}$, $j=1,2,3$.

Recall that if some $\alpha_i$ equals 1 then the corresponding fiber is not singular and the manifold $(b, (O_1,0), (\alpha_1,\beta_1), (\alpha_2,\beta_2), (\alpha_3,\beta_3))$ is not small Seifert manifold but a lens space. Here lens space means a manifold  of genus one or a three-sphere. Well known equivalences for such cases are the following.

A manifold $(b, (O_1,0), (\alpha_1,\beta_1), (\alpha_2,\beta_2), (1,\beta_3))$ is homeomorphic to the manifold $(b~+~\beta_3, (O_1,0), (\alpha_1,\beta_1), (\alpha_2,\beta_2))$.

A manifold $(0, (O_1,0), (\alpha_1,\beta_1), (\alpha_2,\beta_2))$ is homeomorphic to the lens space $L_{p,q}$ with $p=\alpha_1\beta_2+\beta_1\alpha_2$,
$q=\alpha_2\lambda-\beta_2\mu$, where $\alpha_1\lambda+\beta_1\mu=1$.

If some $\alpha_i$ equals 0 then  the presentation $(b, (O_1,0), (\alpha_1,\beta_1), (\alpha_2,\beta_2), (\alpha_3,\beta_3))$ does not present a Seifert fibration but a connected sum of two lens spaces. If $(\alpha_3,\beta_3)=(0,1)$ then the manifold is homeomorphic to
$L_{\alpha_1,\beta_1}\sharp L_{\alpha_2,\beta_2}$.

Consider a such product of three vectors. Let's vary first vector and leave the second and the third vectors fixed.  Then we have the 12 corresponding values of the $t$-invariant. Denote them $t_{\alpha,\beta}$ where ${\alpha,\beta}$ are parameters of the first vector. From the linearity it follows that if a linear combination of some vectors equals null-vector then the same linear combination of the corresponding products equals zero too. So we have an equation for the corresponding values of the $t$-invariant.

There are only two cases of collinear vectors formed by vertices of the icosahedron. The first one is the case of opposite edges. The second one is the case of trapezium formed by four vertices.
The first case implies $t_{\alpha,\beta}+t_{2\alpha,2\beta} = const$ for any $\alpha,\beta$. The second case implies
$t_{\alpha_1,\beta_1}-t_{\alpha_2,\beta_2} = \varepsilon(t_{\alpha_3,\beta_3}-t_{\alpha_4,\beta_4})$ where the vertices $V_1,V_2$ form long base of the trapezium, the  vertices $V_3, V_4$ form short base of the trapezium, and  the vector $V_3V_4$ have the same direction as the vector $V_1V_2$.

Note that $t(S^3)=1$ since $t(M_1 \sharp M_2)=t(M_1)t(M_2)$ ~([2]), $t(RP^3)=\varepsilon+1$ easily follows from the $t$-invariant definition since $RP^2$ is simple spine of $RP^3$.

Consider $t_{\alpha,\beta} = t((-1, (O_1,0), (\alpha,\beta), (0,1), (1,1)))=t(L_{\alpha,\beta})$. We have, that $t(L_{0,1})-t(L_{1,1})=\varepsilon(t(L_{2,1})-t(L_{-1,1}))$. Hence $t(S^2 \times S^1) = \varepsilon +2$.

Also we have $t_{0,1}+t_{0,2} = t_{1,0}+t_{2,0}$. Hence $t(L_{5,2}) = 0$.

If $t_{\alpha,\beta}= t((-1, (O_1,0), (\alpha,\beta), (2,1), (2,1)))$, then $t_{1,1}=t(S^3)=1$, $t_{-1,1}=t(S^3)=1$, $t_{0,1}=t(L_{2,1})t(L_{2,-1})=(\varepsilon+1)^2=3\varepsilon+2$  and $t_{2,1}-t_{-1,1}=\varepsilon^{-1}(t_{0,1}-t_{1,1})$. Therefore $t_{2,1}=\varepsilon+3=t_{2,-1}$.

Also $t_{1,2}+t_{2,-1}=t_{1,0}+t_{2,0}=t_{1,1}+t_{2,2}=t_{0,1}+t_{0,2}$, $t_{-1,2}=\varepsilon+1$,
$t_{1,0}=\varepsilon+2$, $t_{1,1}=1$, $t_{0,1}=3\varepsilon+2$.  Therefore $t_{2,0}=\varepsilon+2$,
$t_{2,2}=2\varepsilon+3$, $t_{0,2}=2-\varepsilon$.

If $t_{\alpha,\beta}= t((-1, (O_1,0), (\alpha,\beta), (2,1), (0,2)))$, then
$t_{1,2}+t_{2,-1}=t_{0,1}+t_{0,2}$, $t_{1,2}=1$, $t_{2,-1}=2-\varepsilon$, $t_{0,1}=0$.
Therefore $t_{0,2}=3-\varepsilon$.

If $t_{\alpha,\beta}= t((-1, (O_1,0), (\alpha,\beta), (0,2), (0,2)))$,
then $t_{1,2}+t_{2,-1}=t_{0,1}+t_{0,2}$, $t_{1,2}=\varepsilon+2$, $t_{2,-1}=3-\varepsilon$, $t_{0,1}=0$.
Therefore $t_{0,2}=5$.

So all cases are considered. Thus the proof is finished.

\vspace{3mm}

{\bf Comparison with the first homology groups}

\vspace{1mm}

$t$-invariant is rather weak already because it has only 12 different values for small Seifert manifolds. However there are some pairs of small Seifert manifolds having the same first homology group but different $t$-invariants. Here we present some examples of that.

\begin{center}
\begin{tabular}{|l|c|c|}  \hline
$b,(\alpha_1,\beta_1),(\alpha_2,\beta_2),(\alpha_3,\beta_3)$
&$H_1$&$t$ \\ \hline
$-1,(2,1),(2,1),(2,1)$&$ Z_2\times Z_2$&$ 3+ \varepsilon$ \\ \hline
$-1,(2,1),(2,1),(4,1)$&$ Z_2\times Z_2$&$       1       $ \\ \hline
$-1,(2,1),(2,1),(3,1)$&$     Z_4      $&$ 3+ \varepsilon$ \\ \hline
$-1,(2,1),(2,1),(5,1)$&$     Z_4      $&$ 2+3\varepsilon$ \\ \hline
$-1,(2,1),(2,1),(3,2)$&$     Z_8      $&$ 3+2\varepsilon$ \\ \hline
$-1,(2,1),(2,1),(5,2)$&$     Z_8      $&$ 2 -\varepsilon$ \\ \hline
\end{tabular}
\end{center}

\begin{center}
\textbf{References}
\end{center}

\lit{1}{Turaev V.\,G., Viro O.\,Y.} {State sum invariants of 3-manifolds and quantum 6j-symbols~// Topology. 1992. V.~31. pp.~865--902.}

\lit{2}{Matveev S.\,V., Ovchinnikov M.\,A., Sokolov M.\,V.} {Construction and properties of the
$t$-invariant~// Journal of Mathematical Sciences. 2003. V. 113. no. 6. pp. 849-855.
[Translated from Zapiski Nauchnykh Seminarov POMI. 2000. V. 267. pp. 207--219.] }

\lit{3}{Ovchinnikov M.\,A.} {Construction of simple spines of Waldhausen manifolds~//
Proceedings of International Conference  "Low-dimensional Topology and Combinatorial Group Theory. Chelyabinsk 1999". Kiev: Institute of Mathematics, 2000. pp.~65--86 (in Russian).}

\lit{4}{Ovchinnikov M.\,A.} {Properties of Viro-Turaev representations of the mapping class group of a torus~// Journal of Mathematical Sciences. 2003. V. 113. no. 6. pp. 856--867. [Translated from Zapiski Nauchnykh Seminarov POMI. 2000. V. 267. pp. 220-240].}

\end{document}